\documentclass[12pt,oneside]{amsart}
\usepackage{amsthm,verbatim,mathtools}
\usepackage{amssymb,latexsym}
\usepackage{amsmath}
\usepackage{graphicx}
\usepackage{color,cite}
\usepackage[dvipsnames]{xcolor}
\usepackage{cancel}
\usepackage{soul}
\usepackage{enumerate}
\newtheorem{theorem}{Theorem}[section]

\newtheorem{remark}[theorem]{Remark}

\newtheorem{corollary}[theorem]{Corollary}
\newtheorem{proposition}[theorem]{Proposition}

\newcommand{\K}{\mathbb K}

\newcommand{\F}{\mathbb F}

\def\daniele#1 {\fbox {\footnote {\ }}\ \footnotetext { From Daniele: {\color{blue}#1}}}

\begin{document}
\title{Permutation trinomials over $\mathbb{F}_{q^3}$}
\author{Daniele Bartoli}
\address{Dipartimento di Matematica e Informatica, Universit\`a degli Studi di Perugia, Via Vanvitelli 1, Perugia, 06123   Italy}

\email {daniele.bartoli@unipg.it}

\maketitle

\begin{abstract}
We consider four classes of polynomials over the fields $\mathbb{F}_{q^3}$, $q=p^h$, $p>3$,  $f_1(x)=x^{q^2+q-1}+Ax^{q^2-q+1}+Bx$, $f_2(x)=x^{q^2+q-1}+Ax^{q^3-q^2+q}+Bx$, 
$f_3(x)=x^{q^2+q-1}+Ax^{q^2}-Bx$,
$f_4(x)=x^{q^2+q-1}+Ax^{q}-Bx$, where $A,B \in \mathbb{F}_q$. We determine conditions on the pairs $(A,B)$ and we  give lower bounds on the number of pairs $(A,B)$ for which these polynomials permute $\mathbb{F}_{q^3}$. 
\end{abstract}


\section{Introduction}

Let $q=p^h$ be a prime power and let $\mathbb{F}_{q}$ denote the finite field with $q$ elements. A polynomial $f(x)\in \mathbb{F}_q[x]$  is said a permutation polynomial (PP for short) if it permutes $\mathbb{F}_q$. Each permutation of the field $\mathbb{F}_q$ can be written as polynomial. On the other hand, in many applications permutation polynomials must satisfy particular constraints or have particular (easy) shapes. 

However, to determine  polynomials belonging to a particular class and permuting the field $\mathbb{F}_q$ can be a difficult task. For a deeper treatment of the connections of PPs with other fields of mathematics we refer to \cite{MuPa, Hou2015} and the references therein. 

One can use different criterions to show that a fixed polynomial permutes a specific $\mathbb{F}_q$. For instance,  Hermite's Criterion states that  a polynomial $f(x)\in \mathbb{F}_q[x]$ is a PP of $\mathbb{F}_q$ if and only $f$ has exactly one root in $\mathbb{F}_q$ and   for each $t\in [1,\ldots,q-2]$, $t \not \equiv 0 \pmod p$, the polynomial $ \left(f(x)\right)^{t} \pmod{x^q-x} $ has degree less than $q-1$; see for instance \cite[Theorem 8.1.7]{MuPa}. In general,  to check the previous condition for every $t$ is quite hard and therefore Hermite's Criterion is mostly used to exclude that a polynomial is a PP. 

Another technique which can be applied to study permutational property of polynomials  relies on the connection between PP and algebraic curves. To a polynomial $f(x) \in \mathbb{F}_{q}[x]$ it can be associated an algebraic curve $\mathcal{C}_f$ of affine equation $(f(x)-f(y))/(x-y)=0$. It is easily seen that $f(x)$ is a PP of $\mathbb{F}_q$ if and only if the curve $\mathcal{C}_f$ does not have points $(a,b)$, $a\neq b \in \mathbb{F}_q$. However, to establish if a particular curve  does not have $\mathbb{F}_{q}$-rational points is a  hard problem. The investigation of the curve $\mathcal{C}_f$ and in particular of its absolutely irreducible components together with the Hasse-Weil Theorem (see \cite[Theorem 5.2.3]{Sti})) give useful criterions and in some cases it can shorten technical proofs; see for instance  \cite{Bartoli2017,BG2017,BMQ2018}. 

Other approaches are applications of the so called AGW criterion; see \cite{AGW2011,YD2011,YD2014}.

Here we use a different technique, based on \cite{Dobbertin2002}, to determine permutation polynomials. Such a method has been applied in \cite{LQC2017,MZFG2017,WZZ2017} in the even characteristic case and in \cite{WZZ2018} for characteristic $3$. 

We study the following four families of permutation polynomials of $\mathbb{F}_{q^3}$ 
$$\begin{array}{lll}
f_1(x)&=&x^{q^2+q-1}+Ax^{q^2-q+1}+Bx,\\
f_2(x)&=&x^{q^2+q-1}+Ax^{q^3-q^2+q}+Bx,\\ 
f_3(x)&=&x^{q^2+q-1}+Ax^{q^2}-Bx,\\
f_4(x)&=&x^{q^2+q-1}+Ax^{q}-Bx,\\
\end{array}
$$
where $A,B\in \mathbb{F}_q$ and the characteristic $p$ is larger than $3$. To prove that a particular polynomial $f_i(x)$ permutes $\mathbb{F}_{q^3}$ we show that for every $a \in \mathbb{F}_{q^3}$ the equation $f_i(x)=a$ has at most one solution in $\mathbb{F}_{q^3}$. We distinguish the cases $a\neq 0$ and $a=0$. 

In the former case we study a system of three equations in $x$, $y=x^q$, and $z=x^{q^2}$ obtained by considering 
\begin{equation}\label{Eq:Intro}
f_i(x)=a, \quad  f_i(x)^q=a^q, \quad f_i(x)^{q^2}=a^{q^2}.
\end{equation}
 In particular, we will show that 
$$\left|\left\{ \overline {x} \in \mathbb{F}_q \ : (\overline{x},\overline{y},\overline{z}) \textrm{ is a solution of } \eqref{Eq:Intro}\right\}\right| \leq 1.$$
This will ensure that  $f_i(x)=a$ has at most one solution.

In the latter case we prove the existence of pairs $(A,B) \in \mathbb{F}_{q}^2$ such that the unique root in $\mathbb{F}_q$ of $f_i(x)=0$ is $x=0$. To this end, we use techniques based on Function Field Theory to provide estimates on the number of $\mathbb{F}_{q}$-rational solutions of particular system of equations. In Section \ref{Sec:AlgGeo} we give a basic introduction on Function Field Theory and we study particular function fields connected with the permutation polynomials $f_i(x)$.

In specific cases, the  computations, mainly related to the resultant of two polynomials, have been done with the help of MAGMA \cite{MAGMA}. The programs are included in the Appendix, whereas an introduction to the use of the resultant can be found in \cite[Section 2]{BG2017}.

\section{Preliminaries from Function Field Theory} \label{Sec:AlgGeo}

In this paper we will make use of some concepts concerning Function Field Theory. This will yield  a lower bound on the number of permutation polynomials of the desired shape. 
 
We recall that a {\em function field} over a perfect field $\mathbb L$ is an extension $\mathbb F$ of $\mathbb L$ such that  $\mathbb F$ is  a finite algebraic extension of $\mathbb L(\alpha)$, with $\alpha$ transcendental over $\mathbb L$. For basic definitions on  function fields we refer to \cite{Sti}. In particular, the (full) constant field of $\mathbb F$ is the set of elements of $\mathbb F$ that are algebraic over $\mathbb L$.

If $\F^{\prime}$ is a finite extension of $\F$, then a place $P^{\prime}$ of $\F^\prime$ is said to be \emph{lying over} a place $P$ of $\F$ if $P \subset P^\prime$. This holds precisely when
$P=P^{\prime} \cap \F$. In this paper, $e(P^{\prime}|P)$ will denote the ramification index of $P^{\prime}$ over $P$. 
A finite extension $\F^{\prime}$ of a function field $\F$ is said to be \emph{unramified} if $e(P^{\prime}|P)=1$ for
every $P^{\prime}$ place of $\F^{\prime}$ and every $P$ place of $\F$ with $P^{\prime}$ lying over $P$. Throughout the paper, we will refer to the following results. 

\begin{theorem}\label{Th_Kummer}\emph{\cite[Cor. 3.7.4]{Sti}}
Consider an algebraic function field $\F$ with constant field $\mathbb L$ containing a primitive $n$-th root of unity ($n>1$ and $n$ relatively prime to the characteristic of $\mathbb L$). Let $u\in \F$ be such that 
there is a place $Q$ of $\F$ with $\gcd(v_Q(u),n)=1$. Let $\F^{\prime}=\F(y)$ with $y^n=u$. Then
\begin{enumerate}
\item $\Phi(T)=T^n-u$ is the minimal polynomial of $y$ over $\F$. The extension $\F^{\prime}:\F$ is Galois of degree $n$ and the Galois group of  $\F^{\prime}:\F$ is cyclic;
\item $$e(P^{\prime}|P)=\frac{n}{r_P} \quad \textrm{where} \quad r_P:= GCD(n,v_P(u))>0\,;$$
\item $\mathbb L$ is the  constant field of $\F^\prime$;
\item let $g^{\prime}$ (resp. $g$) be the genus of $\F^{\prime}$ (resp. $\F$), then
$$ g^{\prime}=1+n(g-1)+\frac{1}{2}\sum_{P\in\mathbb{P}(\F)}(n-r_P)\deg P\,. $$
\end{enumerate}
\end{theorem}
An extension such as $\F^{\prime}$ in Theorem \ref{Th_Kummer}  is said to be a Kummer extension of $\F$.

Denote by $\F_q$ the finite field with $q$ elements and let $\K$ be the algebraic closure of  $\F_q$. A curve $\mathcal{C}$ in some affine or projective space over $\K$ is said to be defined over $\mathbb{F}_q$ if the ideal of $\mathcal{C}$ is generated by polynomials with coefficients in $\mathbb{F}_q$. Let $\K(\mathcal{C})$ denote the function field of $\mathcal C$. The subfield  $\K(\mathcal{C})$ of  $\K(\mathcal{C})$ consists of the rational functions on $\mathcal C$ defined over $\mathbb{F}_q$. The extension $\K(\mathcal C):\F_q(\mathcal C)$ is a constant field extension (see \cite[Section 3.6]{Sti}). In particular, $\mathbb{F}_q$-rational places of $\F_q(\mathcal C)$ can be viewed as the restrictions to $\F_q(\mathcal C)$ of  places of $\K(\mathcal{C})$ that are fixed by the Frobenius map on $\K(\mathcal{C})$. The center of an $\mathbb{F}_q$-rational place is an $\mathbb{F}_q$-rational point of $\mathcal{C}$; conversely,  if $P$ is a simple $\mathbb{F}_q$-rational point of $\mathcal{C}$, then the only place centered at $P$ is $\mathbb{F}_q$-rational. Through the paper, we sometimes use concepts from both Function Field Theory and Algebraic Curves. Concepts as the valuation of a function at a place can be also seen as multiplicity of intersections of fixed algebraic curves; see \cite{Sti}. 

We now recall the well-known Hasse-Weil bound.

\begin{theorem} \label{HasseWeil}\emph{(Hasse-Weil bound, \cite[Theorem 5.2.3]{Sti})}
The number $N_q$ of $\mathbb{F}_q$-rational places of a function field $\,\F$ with constant field $\F_q$ and genus $g$ satisfies 

$$|N_q - (q + 1)| \leq 2g\sqrt{q}.$$
\end{theorem}

\subsection{Some particular function fields associated with permutation polynomials}
In what follows we study particular function fields to obtain lower bounds on the number of the permutation polynomials presented in the subsequent sections. 
\begin{proposition}\label{prop1}
Let $\mathbb{F}_1=\K(x,y)$ be the function field defined by   $x^2 + 3 y^2 - 1=0$, and  $\mathbb{F}=\K(x,y,z,u)$, $q=p^h$, $p>3$, be the function field defined by  the equations
\begin{equation}\label{Eq:EqProp1}
\left\{
\begin{array}{l}
u^3+z^2-z+1=0\\
(6 x^2 y - 6 y^3)z^2+(3 x^3 + 3x^2y - 27 xy^2 - 3y^3 + 3)z -12 x^2y + 12 y^3=0\\
x^2 + 3 y^2 - 1=0\\
\end{array}
\right..
\end{equation}
Then the constant field of $\mathbb{F}$ is $\mathbb{F}_q$ and its genus is at most $11$.  Also, the number of $\mathbb{F}_q$-rational places of $\mathbb{F}$ lying on places of $\mathbb{F}_1$ which are totally unramified is at least
$$q-22\sqrt{q}+1-6\cdot 3 -14=q-22\sqrt{q}-31.$$
\end{proposition}
\proof
First of all observe that Equations \eqref{Eq:EqProp1} are equivalent to 
\begin{equation*}
\left\{
\begin{array}{l}
u^3+z^2-z+1=0\\
\left(z+ \frac{3 x^3 + 3x^2y - 27 xy^2 - 3y^3 + 3}{12(x^2y-y^3)}\right)^2=\frac{\Delta}{36(x^2y-y^3)^2}\\
x^2 + 3 y^2 - 1=0\\
\end{array}
\right.,
\end{equation*}
where $$\Delta=9 x^6 + 18 x^5 y + 135 x^4 y^2 - 180 x^3 y^3 + 18 x^3 + 
    135 x^2 y^4 + 18 x^2 y + 162 x y^5 - 162 x y^2 + 
    297 y^6 - 18 y^3 + 9,$$
that is, recalling that $x^2=-3y^2+1$,
\begin{equation}\label{Eq:Curva2Ter}
\left\{
\begin{array}{l}
u^3+v^2+\frac{12 x y^2 - x + 12 y^3 - 3 y - 1}{2(y-4y^3)}v-\left(\frac{12 x y^2 - x + 4 y^3 - y - 1}{4(y-4y^3)}\right)^2-\frac{12 x y^2 - x + 4 y^3 - y - 1}{4(y-4y^3)}+1=0\\
v^2=\frac{48 x y^5 - 16 x y^3 - 12 x y^2 + x y + x + 48 y^6 - 24 y^4 - 4 y^3 + 3 y^2 
        + y + 1}{2(y-4y^3)^2}\\
x^2 + 3 y^2 - 1=0\\
\end{array}
\right..
\end{equation}
In the following we will consider the function field $\K(x,y,u,v)$ defined by Equations \eqref{Eq:Curva2Ter}.

Consider the function field $\mathbb{F}_1=K(x,y)$ defined by  $x^2 + 3 y^2 - 1=0$. It is a rational function field, since the equation $x^2 + 3 y^2 - 1=0$ defines an absolutely irreducible conic. The function field $\mathbb{F}_2=\mathbb{F}_1(v)$ is an extension of $\mathbb{F}_1$. We want to apply Theorem \ref{Th_Kummer} to show that $\mathbb{F}_2 : \mathbb{F}_1$ is a Kummer extension of degree $2$. To this end, we need to show that  $\Delta$ is not a square in $\mathbb{F}_1$. It is enough to prove the existence of a place $P$ of $\mathbb{F}_1$ such that $v_{P}(\Delta)$ is odd. The resultant $R(y)$ between $\Delta$ and $x^2+3y^2-1$ with respect to $x$ is 
$$y^2(2y-1)^2(2y+1)^2(192 y^6 - 96 y^4 - 80 y^3 + 12 y^2 + 20 y + 11).$$
Recall that $\sum_{P}v_{P}(\Delta)$, $P$ a place centered at $(x_0,y_0)$, corresponds to the multiplicity of intersection between the curves of affine equation $x^2 + 3 y^2 - 1=0$ and $\Delta=0$ at $(x_0,y_0)$.

First of all, we note that the polynomial $192 y^6 - 96 y^4 - 80 y^3 + 12 y^2 + 20 y + 11$ is not a square. In fact, suppose that 
$$192 y^6 - 96 y^4 - 80 y^3 + 12 y^2 + 20 y + 11=192(y^3+ay^2+by+c)^2.$$
Then $a=0$, $b=-1/4$, and both $c=-5/24$ and $192c^2 = 11$, which is impossible since $p>3$.

The existence of a root $\overline{y}$ of odd multiplicity of $R(y)$ yields the existence of at least one place centered at  $(\overline{x},\overline{y})$, for some $\overline{x}$, such that the $v_{P}(\Delta)$ is odd. Therefore, by Theorem \ref{Th_Kummer}, $\mathbb{F}_2$ is a Kummer extension of $\mathbb{F}_1$ and the fields of constants is $\mathbb{F}_q$. The genus of $\mathbb{F}_2$ is 
$$ g^{\prime}=1+2(0-1)+\frac{1}{2}\sum_{P\in\mathbb{P}(\F)}(n-r_P)\deg P\,\leq 2, $$
since the places $P$ for which $v_{P}(\Delta)$ is odd are at most $6$. Consider the place $P$ centered at $(1,0)$ in $\mathbb{F}_1$. Since the valuation $v_P(\Delta)=-2$ there are precisely two places lying over it in $\mathbb{F}_2$. Let $\overline{P}$ one of them. Also, $y$ is a local parameter at $\overline{P}$ and $v_{\overline{P}}(v)=-1$. Since the expansion of $v$ at $\overline{P}$ is $\frac{1}{2y}+\frac{1}{2}+\ldots$ it  can be easily seen that 

$$\alpha=v^2+\frac{12 x y^2 - x + 12 y^3 - 3 y - 1}{2(y-4y^3)}v-\left(\frac{12 x y^2 - x + 4 y^3 - y - 1}{4(y-4y^3)}\right)^2-\frac{12 x y^2 - x + 4 y^3 - y - 1}{4(y-4y^3)}+1$$
$$=\frac{1}{4y^2}-\frac{1}{2y}\cdot \frac{1}{2y}-\frac{1}{16y^2}+\delta,$$
where $v_{\overline{P}}(\delta)>-2$. This means that $v_{\overline{P}}(u^3)=-2$ and therefore by Theorem \ref{Th_Kummer} the extension $\mathbb{F}_3=\mathbb{F}_2(u) : \mathbb{F}_2$ defined by 
$$u^3+v^2+\frac{12 x y^2 - x + 12 y^3 - 3 y - 1}{2(y-4y^3)}v-\left(\frac{12 x y^2 - x + 4 y^3 - y - 1}{4(y-4y^3)}\right)^2-\frac{12 x y^2 - x + 4 y^3 - y - 1}{4(y-4y^3)}+1=0$$
is a Kummer extension. Also, the poles of $v$ in $\mathbb{F}_2$ are at most $3$, whereas the poles of $x$ and $y$ in $\mathbb{F}_2$ are at most $4$ and therefore the poles (and the zeros) of $\alpha$ are at most $7$. Summing up, the genus of $\mathbb{F}_3$ is  
$$ g^{\prime\prime}=1+3(2-1)+\frac{1}{2}\sum_{P\in\mathbb{P}(\F)}(n-r_P)\deg P\,\leq 11. $$
Finally, observe that the number of $\mathbb{F}_q$-rational places of $\mathbb{F}_3$ lying on places of $\mathbb{F}_1$ which are totally unramified is at least
$$q-22\sqrt{q}+1-6\cdot 3 -14=q-22\sqrt{q}-31.$$
\endproof

\begin{remark}\label{Remark}
The number of $\mathbb{F}_q$-rational places of $\mathbb{F}$ lying on places of $\mathbb{F}_1$ which are totally unramified is related to the ``real" number of solutions of System \eqref{Eq:EqProp1}. In fact, we want a lower bound on the number of distinct  centers $(\overline{x},\overline{y},\overline{z},\overline{u})\in \mathbb{F}_q^4$ of the $\mathbb{F}_q$-rational places of $\mathbb{F}$. To this end, we compute only the number of places in $\mathbb{F}$ which are $\mathbb{F}_q$-rational and with ramification index $1$ in the extension $\mathbb{F} : \mathbb{F}_1$.

\end{remark}

\begin{proposition}\label{prop2}
Let $\mathbb{F}_1=\K(x,y)$ be defined by $x^3+y^2-y+1=0$ and  $\mathbb{F}=\K(x,y,z)$, $q=p^h$, $p>3$, $q\equiv 1 \pmod 3$, be the function field defined by  the equations
\begin{equation*}
\left\{
\begin{array}{l}
x^3+y^2-y+1=0\\
\\
z^3=\frac{2y^2 - 3y\alpha + y - 4}{2y^2 + 3y\alpha + y - 4}\\

\end{array}
\right.,
\end{equation*}
where $\alpha^2=-3$. Then the constant field of $\mathbb{F}$ is $\mathbb{F}_q$ and its genus is $2$. Also, there are at least $q-6\sqrt{q}-1$
places of $\mathbb{F}$ lying on places of $\mathbb{F}_1$ which are totally unramified.
\end{proposition}
\proof
It is easily seen that the function field $\mathbb{F}_1=\K(x,y)$ defined by $x^3+y^2-y+1=0$ has genus $1$. The function field $\mathbb{F}$ is therefore a Kummer extension of $\mathbb{F}_1$ since 
$$\gamma=\frac{2y^2 - 3y\alpha + y - 4}{2y^2 + 3y\alpha + y - 4}$$
is not a cube in $\mathbb{F}_1$. To see this it is enough to observe that $\gamma$ has one pole and one zero which are not zeros of $y^2-y+1$. This means that the valuation of $\gamma$ on them is $\pm1$ and therefore $\gamma$ is not a cube in $\mathbb{F}_1$. By Theorem \ref{Th_Kummer}, the constant field of $\mathbb{F}$ is $\mathbb{F}_q$ and its genus is
$$g^{\prime}=1+3(1-1)+ 2=3.$$
Finally, the number of places $P$ in $\mathbb{F}$ lying on places of $\mathbb{F}_1$ which are not totally unramified is at most $2$. This gives a lower bound of 
$$q-6\sqrt{q}-1$$
places of $\mathbb{F}$ lying on places of $\mathbb{F}_1$ which are totally unramified.
\endproof

\begin{proposition}\label{prop3}
Let $q\equiv 1 \pmod 3$. The function field $\mathbb{F}=\K(x,y)$ defined by  
$$Bx^{12} + 12Bx^9 + 24x^9 - 162Bx^6 - 324Bx^3 - 648x^3 + 729B = By^3,$$
with $B^2+B+1=0$ has constant field $\mathbb{F}_q$ and genus at most $4$. Also, there are at least $q-8\sqrt{q}+1-12=q-8\sqrt{q}-11$
places of $\mathbb{F}$ lying over places $P$ of $\K(x)$ which are totally unramified. 
\end{proposition}
\proof
The function field extension $\K(x,y) : \K(x)$ is a Kummer extension if and only if $Bx^{12} + 12Bx^9 + 24x^9 - 162Bx^6 - 324Bx^3 - 648x^3 + 729B$ is not a cube in $\K(x)$. It is easily seen that such a polynomial in $x$ is never a cube. Therefore, by Theorem \ref{Th_Kummer}, the constant field is $\mathbb{F}_q$ and its genus is at most 
$$ g^{\prime}=1+3(0-1)+\frac{1}{2}\sum_{P\in\mathbb{P}(\F)}(n-r_P)\deg P\,\leq 4. $$
Also, there are at least $$q-8\sqrt{q}+1-12=q-8\sqrt{q}-11$$
places $P^{\prime}$ in $\mathbb{F}$ lying over places $P$ of $\K(x)$ which are totally unramified. 

\endproof

\section{Permutation trinomials of type $f_i(x)$, $i=1,2,3,4$}
In this section we investigate permutation trinomials over $\mathbb{F}_{q^3}$ of the form $$\begin{array}{lll}
f_1(x)&=&x^{q^2+q-1}+Ax^{q^2-q+1}+Bx,\\
f_2(x)&=&x^{q^2+q-1}+Ax^{q^3-q^2+q}+Bx,\\ 
f_3(x)&=&x^{q^2+q-1}+Ax^{q^2}-Bx,\\
f_4(x)&=&x^{q^2+q-1}+Ax^{q}-Bx,\\
\end{array}
$$ where $A,B \in \mathbb{F}_q$. 

In what follows, we denote by $\mu_{q^2+q+1}$ the set of the $(q^2+q+1)$-roots of unity in $\mathbb{F}_{q^3}$. 
\begin{theorem}\label{Th:Family1}
Let $A,B \in \mathbb{F}_q$ be such that $A^3+B^2-B+1=0$, $B\neq 0,1$. Suppose that $T^3 + A^2T^2 + (AB + A)T - 1\in \mathbb{F}_{q}[T]$ has no roots in $\mu_{q^2+q+1}$. Then the polynomial  
$$f_{1}(x)=x^{q^2+q-1}+Ax^{q^2-q+1}+Bx$$
is a permutation of $\mathbb{F}_{q^3}$.
\end{theorem}
\proof
Consider $a \in \mathbb{F}_{q^3}$. We will show that the equation 
\begin{equation}\label{Eq1}
x^{q^2+q-1}+Ax^{q^2-q+1}+Bx=a
\end{equation}
has at most one solution. If $a=0$ then $x(x^{q^2+q-2}+Ax^{q^2-q}+B)=0$. If $x\neq 0$ then  
$$u^{q+2}+Au^q+B=0,$$
where $u=x^{q-1}\in \mu_{q^2+q+1}$. Then $u^q=-\frac{B}{u^2+A}$ and $u^{q^2}=\frac{-B(u^2+A)^2}{A(u^2+A)^2+B^2}$. 

Note that $u^2= -A$ would imply $B=0$, a contradiction. Also, $A(u^2+A)^2+B^2=0$ yields $u^{2q}=-A$, that is  $u^2= -A$ again. 

By $u^{q^2+q+1}=1$ and recalling that $A^3+B^2-B+1=0$, we obtain 
$$g_{A,B}(u)=-A u^4 + B^2 u^3 -2 A^2 u^2 +A B^2 u - B + 1=0.$$
To the power $q$ it gives $h_{A,B}(u)=0$, where
$$
\begin{array}{lll}
h_{A,B}(u)&=&(-B + 1)u^8+(-A B^3 - 4 A B + 4 A)u^6+(-3 A^2 B^3 - 2 A^2 B^2 - 6 A^2 B + 6 A^2)u^4\\
&&(-3 A^3 B^3 - 4 A^3 B^2 - 4 A^3 B + 4 A^3 - B^5)u^2-A(A^3 B^3 + 2 A^3 B^2 + A^3 B - A^3 + B^5 + B^4).\\
\end{array}
$$
Now,
$$A^5h_{A,B}(u)+( -A^4 (B-1) u^4-A^3 (B^3-B^2) u^3+A^2 (-A^3 B^3 - 2 A^3 B + 2 A^3 - B^5 + B^4) u^2+$$
$$A(-A^3 B^5 - A^3 B^3 + A^3 B^2 - B^7 + B^6)u-A^6 B^3 - 2 A^6 B^2 - 2 A^6 B + 2 A^6 - A^3 B^7 + A^3 B^2 - 2 A^3 B + A^3$$
$$ - B^9 + B^8)g_{A,B}(u)=B^6(u^3+A^2 u^2 + (A B +A)u - 1).$$
By our assumption there are no $u\in \mu_{q^2+q+1}$ such that  $u^3+A^2 u^2 + (A B +A)u - 1=0$ and therefore $x=0$ is the unique solution of $x(x^{q^2+q-2}+Ax^{q^2-q}+B)=0$.

Suppose now $a\neq 0$. Let $y=x^q$ and $z=y^q$. A solution of Equation \eqref{Eq1} satisfies
\begin{equation}\label{Eq2}
\left\{
\begin{array}{l}
y^2z+A x^2 z+B x^2 y-a x y=0\\
Axy^2 + xz^2 + By^2z - a^qy z=0\\
x^2 y + Bxz^2 - a^{q^2}x z + Ayz^2=0\\
\end{array}
\right.
\end{equation} 
After easy computations (see Appendix), it is easily seen that, for certain $\alpha=\alpha(a,A,B)$ and $\beta=\beta(a,A,B)$,   
$$x^4(Bx - a)(\alpha x+\beta)$$
must vanish. The solution $x=0$ is not possible since $a\neq0$. The solution $x=a/B$ would give $x^{q^2+q-1}+Ax^{q^2-q+1}=0$ and therefore $A^3+1=0$, a contradiction. This means that Equation \eqref{Eq1} has at most one solution.  
\endproof

\begin{corollary}\label{CorFamiglia1}
Let $q=p^h$, $p>3$. Then there are at least $\frac{q-22\sqrt{q}-79}{6}$ pairs $(A,B)\in \mathbb{F}_q^2$ such that $f_{1}(x)$ is a permutation of $\mathbb{F}_{q^3}$.
\end{corollary}
\proof
By Theorem \ref{Th:Family1} it is enough to show prove that there are at least $\frac{q-6\sqrt{q}-56}{3}$ pairs $(A,B)\in \mathbb{F}_q^2$ such that $y^3 + A^2T^2 + (AB + A)T - 1$ has no roots in $\mu_{q^2+q+1}$. 
\begin{itemize}
\item Let $q\equiv 1 \pmod 3$. In this case $-3$ is a square in $\mathbb{F}_q$. By \cite[Theorem 1.34]{HrsBook} the polynomial 
\begin{equation}\label{PolinomioF}
F(T)=y^3 + A^2T^2 + (AB + A)T - 1
\end{equation}
 has three roots in $\mathbb{F}_q$ if and only if the two  roots $\beta_1,\beta_2$ of the Hessian are in $\mathbb{F}_q$ and $F(\beta_1)/F(\beta_2)$ is a cube in $\mathbb{F}_q$. The Hessian of $F(T)$ is 
$$H(T) =(A B^2 + 2A B + 4A)T^2 + (B^3 - 8)T - A^2B^2 - 2A^2B - 4A^2$$
and its roots are
$$t_{1,2}=\frac{-(B^3 - 8)\pm \alpha B(B^2 + 2 B + 4)}{2(A B^2 + 2A B + 4A)},$$
where $\alpha^2=-3$ and $\alpha \in \mathbb{F}_q$. We have that 
$$\frac{F(\beta_1)}{F(\beta_2)}=\frac{2B^2 - 3B\alpha + B - 4}{2B^2 + 3B\alpha + B - 4}.$$

We are interested in the pairs $(A,B)\in \mathbb{F}_q^2$ such that there exists $\xi \in \mathbb{F}_q$ satisfying

\begin{equation}
\left\{
\begin{array}{l}
A^3+B^2-B+1=0\\
\\
\xi^3=\frac{2B^2 - 3B\alpha + B - 4}{2B^2 + 3B\alpha + B - 4}\\
\end{array}
\right..
\end{equation}

By  Proposition \ref{prop2} there are at least 
$$q-4\sqrt{q}-1-15=q-6\sqrt{q}-16$$
triples $(A,B,\xi)\in \mathbb{F}_q^3$ satisfying the previous constraints. The above bound is obtained taking into account that the  degree of the curve defined by the above equations is $15$ and therefore there are at most $15$ places centered at ideal points, that is $15$ places in the corresponding function field considered in Proposition \ref{prop2} which do not correspond to triples     $(A,B,\xi)\in \mathbb{F}_q^3$; see also Remark \ref{Remark}.  Such $q-6\sqrt{q}-16$ triples correspond to at least $\frac{q-6\sqrt{q}-16}{3}$ pairs $(A,B)$: all such pairs give rise to polynomials $F(T)$ having three roots in $\mathbb{F}_q$; see \cite[Theorem 1.34]{HrsBook}. Now we have to exclude cases for which one of these three roots is still in $\mu_{q^2+q+1}$. Note that $\mu_{q^2+q+1}\cap \mathbb{F}_q=\{a_1,a_2,a_3\}$; for each $a_i$ there are at most $6$ pairs $(A,B)$ such that the corresponding polynomial $F(T)$ has $a_i$ as root. Then we have to exclude at most $18$ pairs. Therefore there are at least $\frac{q-4\sqrt{q}-60}{3}$ pairs $(A,B)$ such that $A^3+B^2-B+1=0$ and the polynomial $F(T)$ has no roots in $\mu_{q^2+q+1}$. The corresponding polynomial $f_{1}(x)$ is a  permutation of $\mathbb{F}_{q^3}$.

\item Let $q\equiv 2 \pmod 3$. In this case $-3$ is not a square in $\mathbb{F}_q$. We can argue as above. The roots $\beta_1,\beta_2$ of the Hessian of the polynomial  $F(T)=y^3 + A^2T^2 + (AB + A)T - 1$ are 
$$\beta_{1,2}=\frac{-(B^3 - 8)\pm \alpha B(B^2 + 2 B + 4)}{2(A B^2 + 2A B + 4A)},$$
where $\alpha^2=-3$ and they belong to $\mathbb{F}_{q^2}$. Since $\alpha \in \mathbb{F}_{q^2}\setminus \mathbb{F}_q$,  $\{1,\alpha\}$ is a basis of $\mathbb{F}_{q^2}$ over $\mathbb{F}_{q}$. The three roots of $F(T)$ belong to $\mathbb{F}_q$ if and only if 
$$\frac{F(\beta_1)}{F(\beta_2)}=\frac{2B^2 - 3B\alpha + B - 4}{2B^2 + 3B\alpha + B - 4} $$
is a cube in $\mathbb{F}_{q^2}$. This means that there exist $\xi_1,\xi_2\in \mathbb{F}_{q}$ such that 
\begin{equation}
\left\{
\begin{array}{l}
A^3+B^2-B+1=0\\
\\
(\xi_1+\alpha\xi_2)^3=\frac{2B^2 - 3B\alpha + B - 4}{2B^2 + 3B\alpha + B - 4}\\
\end{array}
\right.
\iff 
\left\{
\begin{array}{l}
A^3+B^2-B+1=0\\
\\
(2 \xi_1^3 - 18 \xi_1 \xi_2^2 - 2)B^2\\
\hspace{1 cm}+(\xi_1^3 - 27 \xi_1^2 \xi_2 - 9 \xi_1 \xi_2^2 + 27 \xi_2^3 - 1)B\\
\hspace{1 cm} -4 \xi_1^3 + 36 \xi_1 \xi_2^2 + 4=0\\
    \\
(6 \xi_1^2 \xi_2 - 6 \xi_2^3)B^2\\
\hspace{1 cm}+(3 \xi_1^3 + 3 \xi_1^2 \xi_2 - 27 \xi_1 \xi_2^2 - 3 \xi_2^3 + 3)B\\
\hspace{1 cm} -12 \xi_1^2 \xi_2 + 12 \xi_2^3=0\\
\end{array}
\right..
\end{equation}
The previous system is equivalent to

\begin{equation*}
\left\{
\begin{array}{l}
A^3+B^2-B+1=0\\
\\
(6 \xi_1^2 \xi_2 - 6 \xi_2^3)B^2\\
\hspace{1 cm}+(3 \xi_1^3 + 3 \xi_1^2 \xi_2 - 27 \xi_1 \xi_2^2 - 3 \xi_2^3 + 3)B\\
\hspace{1 cm} -12 \xi_1^2 \xi_2 + 12 \xi_2^3=0\\
\\
(\xi_1^2 + 3 \xi_2^2 - 1)(\xi_1^4 + 6 \xi_1^2 \xi_2^2 + \xi_1^2 + 9 \xi_2^4 + 3\xi_2^2 + 1)=0\\
    \\
\end{array}
\right..
\end{equation*}

By  Proposition \ref{prop1}, there are at least 
$$q-22\sqrt{q}-31-30=q-22\sqrt{q}-61$$
 quadruples $(A,B,\xi_1,\xi_2)\in \mathbb{F}^4$ satisfying the previous set of constraints (here we used that the degree of the curve corresponding to the function field of Proposition \ref{prop1} is at most $30$); see also Remark \ref{Remark}. Therefore there are at least 
 $$\frac{q-22\sqrt{q}-61}{6}$$
pairs $(A,B)$ such that $A^3+B^2-B+1=0$ and the polynomial $F(T)$ has  three roots in $\mathbb{F}_q$; see \cite[Theorem 1.34]{HrsBook}. Since $\mu_{q^2+q+1}\cap \mathbb{F}_q=\{1\}$ and $1$ is a root of $F(T)$ if and only if $A^2+AB+A=0$, easy computations show that this happens for at most three pairs $(A,B)$ with $A^3+B^2-B+1=0$. This yields the existence of at least 
$$\frac{q-22\sqrt{q}-79}{6}$$ 
pairs $(A,B)\in \mathbb{F}_q^2$ for which the corresponding polynomial $f_{1}(x)$  permutes  $\mathbb{F}_{q^3}$. 
\end{itemize}
\endproof

\begin{theorem}\label{Family4}
Let $A,B \in \mathbb{F}_q$ be such that $A^3+B^2-B+1=0$, and $B\neq 0,-1$. Then there are at least $\frac{q-22\sqrt{q}-79}{6}$ pairs $(A,B)\in \mathbb{F}_q^2$ such that 
$$f_{2}(x)=x^{q^2+q-1}+Ax^{q^3-q^2+q}+Bx$$
is a permutation of $\mathbb{F}_{q^3}$.
\end{theorem}
\proof

We need to prove that for every $a\in \mathbb{F}_{q^3}$ there exists at most one solution of 
\begin{equation}\label{Eq1Bis}
x^{q^2+q-1}+Ax^{q^3-q^2+q}+Bx=a.
\end{equation}

Suppose  $a\neq 0$. Let $y=x^q$ and $z=y^q$. A solution of the above equation  satisfies
$$
\left\{
\begin{array}{l}
Ax^2y +B x^2 z B - a x z + y z^2=0\\
x^2 z + Bx y^2 - a^qx y + Ay^2 z=0\\
x y^2 + Ax z^2+ By z^2 - a^{q^2}y z=0\\
\end{array}
\right..
$$
After easy computations (see Appendix), it is easily seen that, for certain $\alpha=\alpha(a,A,B)$ and $\beta=\beta(a,A,B)$,
$$x^4(Bx - a)(\alpha x+\beta)$$
must vanish. The solution $x=0$ is not possible since $a\neq0$. The solution $x=a/B$ would give $x^{q^2+q-1}+Ax^{q^3-q^2+q}=0$ and therefore $A^3+1=0$, a contradiction. This means that Equation \eqref{Eq1Bis} has at most one solution if $a\neq 0$.

Suppose now $a=0$. Then, apart from $x=0$, another solution of the previous equation satisfies
$$x^{2q^2+q}+Ax^{q+2}-Bx^{q^2+2}=0 \iff x^{q^2+2}+Ax^{q^2+2q}+Bx^{2q+1}=0 \iff x^{q^2-q}+Ax^{q^2+q-2}+Bx^{q-1}=0.$$
Let us call $y=x^{q-1}$, then $ y^{q}+Ay^{q+2}+By=0$, that is 
$$y^q=\frac{-By}{Ay^2+1}, \qquad y^{q^2}=\frac{B^2y(Ay^2+1)}{AB^2y^2+(Ay^2+1)^2}.$$

Note that $Ay^2=-1$ would yield $B=0$, a contradiction. On the other hand, if $AB^2y^2+(Ay^2+1)^2=0$ then, by $ y^{q}+Ay^{q+2}+By=0$, we obtain $Ay^{2q}+1=0$, that is $Ay^{2}+1=0$ again. 

Since $y^{q^2+q+1}=1$ we get, using $A^3+B^2-B+1=0$,
$$M(y)=-A^2y^4 -B^3y^3 -AB^2 y^2 - 2Ay^2 - 1=0.$$
Rising to the power $q$ and using $y^q=\frac{-By}{Ay^2+1}$ we obtain 

$$
\begin{array}{rrr}
L(y)=-A^4y^8 + (-A^3B^4 - 2A^3B^2 - 4A^3)y^6+AB^6y^5\\
     +(- 3A^2B^4 - 4A^2B^2 - 6A^2)y^4 +B^6y^3+(- AB^4 - 2AB^2 - 4A)y^2 - 1&=&0.
\end{array}
$$
Since 
$$
\begin{array}{l}
A^6 L(y) +\Big(-A^8 y^4+A^6 B^3 y^3+A^4 (-A^3 B^4 - A^3 B^2 - 2 A^3 - B^6) y^2\\
-A^2 B^6 (-A^3 B - A^3 - B^3) y+(A^3 B^2 - A^3 B + A^3 - B^6) (A^3 B^4 + A^3 B^3 - A^3 B - A^3 + B^6)\Big) M(y)\\
=B^6(y^3 - ABy^2 - Ay^2 - A^2y - 1),\\
\end{array}
$$
it is enough to show that the polynomial $G(T)=T^3 - ABT^2 - AT^2 - A^2T - 1$ has no roots in $\mu_{q^2+q+1}$. Note that the polynomial $-T^3G(1/T)$ is exactly the polynomial $F(T)$ in \eqref{PolinomioF}. Therefore arguing as in Corollary \ref{CorFamiglia1} the assertion follows.

\endproof

Let 

\begin{equation}\label{h_1}
\begin{array}{lll}
h_1(T)&=&T(T^{22} + 9 T^{19} + 18 T^{18} + T^{17} + 63 T^{16} + 78 T^{15} + 72 T^{14} + 165 T^{13} - 215 T^{12} \\
&& - 64 T^{11}+ 300 T^{10} - 108 T^9 + 15 T^8 + 45 T^7 + 7 T^6 + 
        36 T^5 + 6 T^4 + 3 T^2 + 1)\cdot\\
&&(T^{44} - T^{42} + 3 T^{41} - 18 T^{40} - 7 T^{39} + 22 T^{38} - 50 T^{37} + 118 T^{36} + 145 T^{35} \\
&& - 254 T^{34}+ 218 T^{33} - 112 T^{32} - 726 T^{31} + 627 T^{30} - 217 T^{29} - 16 T^{28} + 258 T^{27} \\
&& + 996 T^{26} - 611 T^{25}- 161 T^{24} - 691 T^{23} - 392 T^{22} + 1189 T^{21} + 252 T^{20} - 645 T^{19} \\
&& + 458 T^{18} - 475 T^{17} - 141 T^{16}+ 237 T^{15} + 72 T^{14} + 298 T^{13} - 327 T^{12} -121 T^{11} \\
&& + 140 T^{10} + 47 T^9 + 27 T^8 - 59 T^7 - 10 T^6+ 22 T^5 + 3 T^4- 3 T^2 + 1).
\end{array}
\end{equation}

\begin{theorem}\label{Family2}
Let $q\equiv 1 \pmod 3$. If $A,B \in \mathbb{F}_q$, $B^2+B+1=0$, $h_1(A)\neq 0$, $A^3\neq -1$, then the polynomial  
$$f_3(x)=x^{q^2+q-1}+Ax^{q^2}-Bx$$
is a PP of $\mathbb{F}_{q^3}$.
\end{theorem}
\proof

Let us denote $x^q$ and $x^{q^2}$ by $y$ and $z$. The polynomial $x^{q^2+q-1}+Ax^{q^2}-Bx$ permutes $\mathbb{F}_{q^3}$ if and only if for every $a \in \mathbb{F}_{q^3}$ there exists at most one solution to $x^{q^2+q-1}+Ax^{q^2}-Bx=a$. 

Let $a=0$. Clearly a solution is $x=0$. Suppose $x\neq0$. Then, putting $y=x^{q-1}$, we have $y^{q+2}+Ay^{q+1}-B=0$, that is 
$$y^q=\frac{B}{y(y+A)},\qquad y^{q^2}=\frac{y^2(y+A)^2}{B+Ay(y+A)}.$$
Note that $y\neq 0$ and $y\neq -A$ otherwise $B=0$. Also, $B+Ay(y+A)=0$ would imply $y(y+A)=-B/A$ and then $y^q=-A$, that is again $y=-A$, a contradiction.
From $y^{q^2+q+1}=1$ we obtain 
$$A^2By^4 - Ay^3 - 2A^2y^2 -A^3 y - By - AB=0$$
and to the power $q$ 
$$
\begin{array}{rrr}
-ABy^8 - 4A^2By^7 - (7A^3B+B^2)y^6  - (7A^4B+3AB^2)y^5 -    4A^5B+5A^2B^2)y^4&&\\ 
- (A^6B+5A^3B^2)y^3-( 2A^4B^2+AB^3)y^2- A^2B^3y + A^2B^5&=&0.\\
    \end{array}$$
    
Recalling that $B^2+B+1=0$ and using the two previous relations on $y$, we obtain that $h_1(A)=0$ a contradiction; see Appendix.

If $a\neq 0$ then $x\neq 0$ and the previous equation can be written as $$yz+Axz-Bx^2-ax=0.$$ Rising to the power $q$ we get $zx+Ayx-By^2-a^qy=0$ and $xy+Azy-Bz^2-a^{q^2}z=0$. It can be seen (see Appendix) that these three conditions imply $x=0$, $x=-a/B$ or 
\begin{center}
${\small x=\frac{(A^2 B  + A^2)a^{q^2+2q+1} + (A B-A) a^{2q^2+q+1} + A B a^{q^2+3q} - A a^{2q+2}  - B a^{q^2+q+2} - (B+1) (a^{3q+1} + a^{3q^2+1}) - a^{2q^2+2q}}{(A^3 B+ A^3- 3 B-3) a^{q^2+q+1} + A^2 B a^{q+2} + A^2 B a^{2q^2+1} + A^2 B a^{q^2+2q} + 
        (A B-A) (a^{q^2+2} +  a^{2q+1} +  a^{2q^2+q})   - a^3  - a^{3q} - a^{3q^2}}}.$
\end{center}
Clearly $x=0$ cannot be a solution. If $x=-a/B$ is a solution then from $x^{q^2+q-1}+Ax^{q^2}-Bx=a$ we get $x^{q^2+q-1}+Ax^{q^2}=0$ which yields $x^{q-1}=-A$ and then $1=x^{q^3-1}=\left(x^{q-1}\right)^{q^2+q+1}=\left(-A\right)^{q^2+q+1}=-A^3$, a contradiction. 

Therefore for each $a\in \mathbb{F}_q$ the equation $x^{q^2+q-1}+Ax^{q^2}-Bx=a$ has at most one solution and the polynomial $f_3(x)$ is a permutation polynomial of $\mathbb{F}_{q^3}$. 
\endproof

\begin{theorem}\label{Family3}
Let $q\equiv 1 \pmod 3$. If $A,B \in \mathbb{F}_q$, $B^2+B+1=0$, $A^3\neq -1$ then there are at least $\frac{q-8\sqrt{q}-50}{3}$ values $A\in \mathbb{F}_q$ for which the  polynomial  
$$f_4(x)=x^{q^2+q-1}+Ax^{q}-Bx$$
is a PP of $\mathbb{F}_{q^3}$.
\end{theorem}
\proof
We have to prove that for every $a \in \mathbb{F}_{q^3}$ there exists at most one solution to $x^{q^2+q-1}+Ax^{q}-Bx=a$. 

Let $a=0$. Clearly a solution is $x=0$. Suppose $x\neq0$. Then, putting $y=x^{q-1}$, we have $y^{q+2}+Ay-B=0$, that is 
$$y^q=\frac{B-Ay}{y^2},\qquad y^{q^2}=\frac{By^4-A(B-Ay)y^2}{(B-Ay)^2}.$$
Clearly $y\neq 0$. If $y=B/A$ then from $y^{q+2}+Ay-B=0$ we get $B=0$, a contradiction. 
From $y^{q^2+q+1}=1$ we obtain 
$$F(y)=By^3 + A^2y^2 - ABy + Ay - B=0.$$

The hessian with respect to $y$ is 
$$H(y)=-y((A^4 - 6AB - 3A)y - A^3B + A^3 - 9B - 9)$$
and its roots are $\beta_1=0$ and $\beta_2=\frac{ A^3B - A^3 + 9B + 9}{A^4 - 6AB - 3A}$. Therefore
$$\frac{F(\beta_1)}{F(\beta_2)}=\frac{B(A^4 - 6AB - 3A)^3}{A^{12} B + 12 A^9 B + 24 A^9 - 162 A^6 B - 324 A^3 B - 648 A^3 + 729 B} $$
is a cube in $\mathbb{F}_{q}$ if and only if 
$$A^{12} B + 12 A^9 B + 24 A^9 - 162 A^6 B - 324 A^3 B - 648 A^3 + 729 B=BC^3$$
for some $C\in \mathbb{F}_{q}$. By Proposition \ref{prop3} there are at least 
$\frac{q-8\sqrt{q}-23}{3}$
elements $A\in \mathbb{F}_q$ such that $\frac{F(\beta_1)}{F(\beta_2)}$ is a cube in $\mathbb{F}_q$ and therefore by \cite[Theorem 1.34]{HrsBook} the three roots of $F(y)$ are in $\mathbb{F}_q$. Since $|\mu_{q^2+q+1}\cap \mathbb{F}_q|=3$ there are at most $6$ values of $A$ for which $F(y)$ has roots in $\mu_{q^2+q+1}\cap \mathbb{F}_q$. This means that for at least 
$$\frac{q-8\sqrt{q}-41}{3}$$
elements $A\in \mathbb{F}_q$ the equation $x^{q^2+q-1}+Ax^{q}-Bx=0$ has $0$ as unique root. 

If $a\neq 0$ then $x\neq 0$ and the previous equation can be written as $$yz+Axz-Bx^2-ax=0.$$ Rising to the power $q$ we get $zx+Ayx-By^2-a^qy=0$ and $xy+Azy-Bz^2-a^{q^2}z=0$. It can be seen (see Appendix) that these three conditions imply $x=0$, $x=-a/B$ or 
\begin{center}

${\small x=\frac{(A^2 B  + A^2)a^{2q^2+q+1} + (A B-A) a^{q^2+2q+1} + A B a^{3q^2+q} - A a^{2q^2+2} - (B+1) (a^{3q+1} + a^{3q^2+1}) - a^{2q^2+2q}-Ba^{q^2+q+2}}{(A^3 B+ A^3- 3 B-3) a^{q^2+q+1} + A^2 B a^{q^2+2} + A^2 B a^{2q+1} + A^2 B a^{2q^2+q} + 
        (A B-A) (a^{q^2+2q} +  a^{q+2} +  a^{2q^2+1})   - a^3  - a^{3q} - a^{3q^2}}}.$
\end{center}
Clearly $x=0$ cannot be a solution. If $x=-a/B$ is a solution then from $x^{q^2+q-1}+Ax^{q}-Bx=a$ we get $x^{q^2+q-1}+Ax^{q}=0$ which yields $x^{q-1}=-A$ and then $1=x^{(q^3-1)(q+1)}=\left(x^{q^2-1}\right)^{q^2+q+1}=\left(-A\right)^{q^2+q+1}=-A^3$, a contradiction. This means that for every $a\in \mathbb{F}_{q^3}$ the equation $x^{q^2+q-1}+Ax^{q^2}-Bx=a$ has at most one solution and $f_4(x)$ is a permutation polynomial of $\mathbb{F}_{q^3}$.

Since there are three values in $\mathbb{F}_q$ which are $3$-rd roots of $-1$ there are at least  
$$\frac{q-8\sqrt{q}-50}{3}$$
elements $A\in \mathbb{F}_q$ for which 
$$x^{q^2+q-1}+Ax^{q}-Bx$$
is a PP of $\mathbb{F}_{q^3}$.
\endproof

\section{Appendix}
We list some MAGMA programs used in the proofs of the results in the previous sections. Note that in general we performed the computations in rings of polynomials over the integers. 

We will make use a number of time of the function ``Substitution"
{\scriptsize
\begin{verbatim}
Substitution := function (pol, m, p)
	e := 0;
	New := K! pol;
	while e eq 0 do
		N := K!0;
		T := Terms(New);
		i:= 0;
		for t in T do
			if IsDivisibleBy(t,m) eq true then
				Q := K! (t/m);
				i := 1;
				N := K!(N + Q* p);
			else 
				N := K!(N + t);
			end if;
		end for;
		if i eq 0 then 
			return New;
		else	
			New := K!N;
		end if;	
	end while;
end function;
\end{verbatim}
}

\subsection{Theorem \ref{Th:Family1}}
{\scriptsize
\begin{verbatim}
K<x,y,z,A,B,a,b,c> := PolynomialRing(Integers(),8);
p1 := y^2*z+A*x^2*z+B*x^2*y-a*x*y;
p2 := Evaluate(p1,[y,z,x,A,B,b,c,a]);
p3 := Evaluate(p2,[y,z,x,A,B,b,c,a]);
R1 := K!(Resultant(p1,p2,z)/x/y^2);
R2 := K!(Resultant(p1,p3,z)/x^2/y);
RR := Resultant(R1,R2,y);
RR := Substitution(RR,A^3,-B^2+B-1);
Factorization(RR);
\end{verbatim}
}

\subsection{Theorem \ref{Family4}}
{\scriptsize
\begin{verbatim}
K<x,y,z,A,B,a,b,c> := PolynomialRing(Integers(),8);
p1 := y*z^2+A*x^2*y+B*x^2*z-a*x*z;
p2 := Evaluate(p1,[y,z,x,A,B,b,c,a]);
p3 := Evaluate(p2,[y,z,x,A,B,b,c,a]);
R1 := K!(Resultant(p1,p2,y)/x^2/z);
R2 := K!(Resultant(p1,p3,y)/x/z^2);
RR := Resultant(R1,R2,z);
RR := Substitution(RR,A^3,-B^2+B-1);
Factorization(RR);
\end{verbatim}
}

\subsection{Theorem \ref{Family2}}
{\scriptsize
\begin{verbatim}
K<x,y,z,A,B,a,b,c> := PolynomialRing(Integers(),8);
p1 := z*y+A*z*x-B*x^2-a*x;
p2 := Evaluate(p1,[y,z,x,A,B,b,c,a]);
p3 := Evaluate(p2,[y,z,x,A,B,b,c,a]);
R1 := K!(Resultant(p1,p2,z));
R2 := K!(Resultant(p1,p3,z)/x);
RR := K!(Resultant(R1,R2,y)/x^2/(x*B + a)^3);
RR := Substitution(RR,B^2,-B-1);
Factorization(RR);
//FOR THE CASE a=0
yq := B/y/(y+A);
yqq := y^2*(y*A)^2/(B+A*y*(y+A));
pol1 := K!((y+A)*(B+A*y*(y+A))*(yqq*yq*y-1));
pol2 := K!((y+A)^4*y^4*Evaluate(pol1,[x, yq,z,A,B,a,b,c]));
pol2 := Substitution(pol2,B^2,-B-1);
RR := Resultant(pol1,pol2,y);
RR := Resultant(RR,B^2+B+1,B);
Factorization(RR);
\end{verbatim}
}

\subsection{Theorem \ref{Family3}}
{\scriptsize
\begin{verbatim}
K<x,y,z,A,B,a,b,c> := PolynomialRing(Integers(),8);
p1 := z*y+A*y*x-B*x^2-a*x;
p2 := Evaluate(p1,[y,z,x,A,B,b,c,a]);
p3 := Evaluate(p2,[y,z,x,A,B,b,c,a]);
R1 := K!(Resultant(p1,p2,z));
R2 := K!(Resultant(p1,p3,z)/x);
RR := K!(Resultant(R1,R2,y)/x^2/(x*B + a)^3);
RR := Substitution(RR,B^2,-B-1);
Factorization(RR);
//FOR THE CASE a=0
yq := B/y/(y+A);
yqq := y^2*(y*A)^2/(B+A*y*(y+A));
pol1 := K!((y+A)*(B+A*y*(y+A))*(yqq*yq*y-1));
pol2 := K!((y+A)^4*y^4*Evaluate(pol1,[x, yq,z,A,B,a,b,c]));
pol2 := Substitution(pol2,B^2,-B-1);
RR := Resultant(pol1,pol2,y);
RR := Resultant(RR,B^2+B+1,B);
Factorization(RR);
\end{verbatim}
}
\section*{Acknowledgments}

The author was partially supported by the Italian Ministero dell'Istruzione, dell'Universit\`a e della Ricerca (MIUR) and by the Gruppo Nazionale per le Strutture Algebriche, Geometriche e le loro Applicazioni (GNSAGA-INdAM).

\end {document}